\documentclass[12pt, letterpaper, reqno]{amsart}
\usepackage[hyperindex,breaklinks]{hyperref}
\usepackage{amsmath}
\usepackage{amsfonts}
\usepackage{mathrsfs}
\usepackage{graphicx}
\usepackage{color}
\usepackage{slashed}
\usepackage{braket}
\usepackage{extarrows}
\usepackage{amsmath, bm}
\usepackage{amssymb}
\usepackage[margin=1.5in]{geometry}
\usepackage{bbm}
\usepackage{hyperref}

\definecolor{yellow}{rgb}{1,.7,0}

\definecolor{pkured}{rgb}{0.55,0,0}

\newcommand{\be}{\begin{equation*}}
	\newcommand{\ee}{\end{equation*}}
\newcommand{\beq}{\begin{equation}}
	\newcommand{\eeq}{\end{equation}}

\numberwithin{equation}{section}

\newtheorem{cor}{Corollary}[section]

\newtheorem{lem}[cor]{Lemma}
\newtheorem{prop}[cor]{Proposition}
\newtheorem{thm}[cor]{Theorem}

\theoremstyle{remark}
\newtheorem{rmk}[cor]{Remark}

\numberwithin{figure}{section}
\usepackage{ytableau}

\usepackage{tikz}
\newcounter{x}
\newcounter{y}
\newcounter{z}

\author{Ce Ji}
\email{cji@mail.tsinghua.edu.cn}
\address{Ce Ji, Department of Mathematical Sciences, Tsinghua University, Beijing, China}

\author{Qian Tang}
\email{919130201@qq.com}
\address{Qian Tang, Department of Mathematics, The University of Hong Kong, Hong Kong SAR, China}

\author{Chenglang Yang}
\email{yangcl@amss.ac.cn}
\address{Chenglang Yang, Hua Loo-Keng Center for Mathematical Sciences,
	Academy of Mathematics and Systems Science,
	Chinese Academy of Sciences,
	Beijing, China}

\title[{Schr{\"o}der Paths, Their Generalizations and Knot Invariants}]
{Schr{\"o}der Paths, Their Generalizations and Knot Invariants}

\begin{document}

\begin{abstract}
    We study some kinds of generalizations of Schr\"oder paths below a line with rational slope and derive the $q$-difference equations that are satisfied by their generating functions.
    As a result,
    we establish a relation between the generating function of generalized Schr\"oder paths with backwards and the wave function corresponding to colored HOMFLY-PT polynomials of torus knot $T_{1,f}$.
    We also give a combinatorial proof of a recent result by Sto\v si\'c and Su\l kowski,
    in which the standard generalized Schr\"oder paths are related to the superpolynomial of reduced colored HOMFLY-PT homology of $T_{1,f}$.
\end{abstract}
\maketitle

\setcounter{section}{0}
\setcounter{tocdepth}{2}


\section{Introduction}

In the study of the knot theory,
various knot invariants, including the Jone polynomials and HOMFLY-PT polynomials, play a central role (see \cite{HOMFLY,Jon87,PT}).
As a variant of the HOMFLY-PT polynomials, colored HOMFLY-PT polynomials were introduced by Witten using Chern-Simons theory \cite{Wit89}, 
and further defined by Reshetikhin and Turaev using quantum groups and R-matrix approach \cite{RT90}.
These knot invariants are closely related to open string theory in Calabi-Yau geometry (see \cite{AKMV,Wit95}),
 which are expected to have counting meanings.

Counting various lattice paths is a widely studied problem in combinatorics,
and is related to many other aspects of mathematics (see, for example, \cite{BF02} and references therein). 
There are many interesting lattice path counting problems, and in this paper we are in particular interested in counting the Schr\"oder paths and their generalizations.
The standard Schr\"oder paths are paths in a square lattice from $(0,0)$ to $(l,l)$ for some integer $l\geq0$ with three candidate steps $(0,1), (1,0)$ and $(1,1)$,
which lie below the line $y=x$.
The generalized Schr{\"o}der paths are considered as paths below the line with rational slope $y=\frac{mx}{n}$.
The counting problems related to Schr\"oder paths are intensively studied.
For examples,
Bizley \cite{Biz54} gave the generating function of numbers of such paths when there is no diagonal step $(1,1)$,
which are called the Dyck paths in literature.
In general,
Schr\"oder paths below a line with rational slope were studied, for example, in \cite{BF02,BW19,Son05} and references therein.
In two recent seminal papers \cite{PSS18,SS24},
Panfil, Sto\v si\'c and Su\l kowski proposed to count the standard generalized Schr{\"o}der paths labeled by certain areas and the number of diagonal steps.
They promised that this counting problem should be related to corresponding torus knot invariants.
Motivated by their paper,
we study these kinds of Schr{\"o}der paths and their generalizations.

We first consider the generalized Schr\"oder paths below the line $y=\frac{mx}{n}$.
We count the generalized Schr{\"o}der paths from $(0,0)$ to $(nl,ml)$ for some integer $l\geq0$, with all of their midway points lying below the line $y=\frac{mx-s}{n}$ for some integers $0\leq s\leq mn$. The generating functions of such countings are defined as $y^{[s]}(a,q;x).$
See Section \ref{subsec:slopemn} for their explicit definitions. 
Our first main result is stated as follows.

\begin{thm}\label{thm:main1}
The series of generating functions $y^{[0]},y^{[1]},\ldots,y^{[mn-1]},y^{[mn]}$ are determined by the following $q$-difference equations
\begin{subequations}
\begin{align}
    y^{[0]}(a,q;x) & = 
    \frac{1}{1-y^{[1]}(a,q;x)},\\
    \nonumber
    y^{[s]}(a,q;x) & =
    y^{[s+1]}(a,q;x) +
    q^{mn\epsilon_s^2-2\alpha_s\beta_s}
    x^{-\epsilon_s}\cdot
    y^{[0]}(a,q;q^{2s}x)\\
    & \quad \quad\ \ \cdot
    y^{[\beta_sn+1]}(a,q;q^{2(s-\beta_sn)}x)\cdot
    y^{[\alpha_sm+1]}(a,q;q^{2(s-\alpha_sm)}x),\\
    y^{[mn]}(a,q;x) & = a^2q^{mn-1}x+
    q^{mn}x\cdot y^{[0]}(a,q;q^{2mn}x)
\end{align}
\end{subequations}
with initial values
\begin{align}
    y^{[0]}(a,q;0)=1,\quad
    y^{[s]}(a,q;0)=0,\  0<s<mn,\quad
    y^{[mn]}(a,q;0)=0,
\end{align}
    where the integers $\alpha_s,\beta_s$ are determined by
    \begin{alignat*}{2}
    0 &\leq \alpha_s < n,&\quad
    \alpha_sm &= s \mod n,\\
    0 &\leq \beta_s < m,&\quad
    \beta_sn &= s \mod m,
    \end{alignat*}
    and $\epsilon_s:=\frac{\alpha_sm+\beta_sn-s}{mn}$.
\end{thm}

When considering the lattice paths under the line $y=\frac{x}{f}$ for positive integer $f$,
we also study other three types of generalized Schr\"oder paths, with corresponding generating functions denoted by 
$$y_{\infty}(a,q;x),\quad \{y_i(a,q;x)\}_{i=1}^{f+1}, \quad h(a,q;x)$$
with following descriptions:
\begin{itemize}
\item $y_i(a,q;x)$ counts generalized Schr{\"o}der paths from $(0,0)$ to $(fl+i-1,l)$ for some integer $l\geq0$;
\item $y_\infty(a,q;x)$ counts generalized Schr{\"o}der paths with finite height and without endpoint,
\item $h(a,q;x)$ counts generalized Schr\"oder paths with backwards,
i.e. this kind of paths admit four candidate paths as $(0,1), (1,0), (1,1)$ and $(-1,0)$.
\end{itemize}
See Section \ref{sec:gSchroder2} for further explicit definitions and properties. 
We can explicitly write down the $q$-difference equations for each of these three types of generating functions.
In particular,
they are related to certain knot invariants for the torus knot $T_{1,f}$ (see \cite{PSS18,SS24}).

For the torus knots 
their colored HOMFLY-PT polynomials are explicitly derived by Rosso and Jones \cite{RJ93} using quantum groups.
Later,
Lin and Zheng \cite{LZ10} gave them a formula in terms of certain combinatorics of Schur polynomials.
The generating function of all these colored HOMFLY-PT polynomials forms a partition function,
which is expected to possess certain integrable systems,
especially the KP integrable hierarchies (see \cite{BEM12, DPSS19, MMS13}).
Inspired by the study of integrable systems, we consider the wave function of the corresponding torus knots $T_{1,f}$,
which can be directly derived from certain specialization of the partition function
and is related to corresponding quantum $A$-polynomial and quantum spectral curve \cite{GS12}.
As a result of our combinatorial study on the Schr\"oder paths and their generalizations,
we first prove the following relation between the generating function of generalized Schr\"oder paths with backwards $h(a,q;x)$ and the wave function $\psi^{T_{1,f}}(x)$ corresponding to the torus knot $T_{1,f}$.

\begin{prop}\label{thm:main wave}
    Denote by $h(a,q;x)$ the generating function of generalized Schr{\"o}der paths with backwards under the line $y=\frac{x}{f}$, 
    and let $\psi^{T_{1,f}}(x)$ be the wave function corresponding to the partition function of colored HOMFLY-PT polynomials of torus knot $T_{1,f}$ (See Section \ref{subsec:wave} for definitions).
    We have
    \begin{align}\label{eqn:main wave}
        h(a,q;x)
        =\psi^{T_{1,f}}(-\nu^{1/2}q^{f-1}x)|_{\nu\rightarrow-a^2/q, t\rightarrow q^{-2}}.
    \end{align}
\end{prop}
We hope the relation in Proposition \ref{thm:main wave} could be extended to the Baker-Akhiezer kernel $\psi(w,z)$,
which is a generalization of the wave function,
since this kernel $\psi(w,z)$ contains all the information of corresponding partition function.
After the change of variables demonstrated in equation \eqref{eqn:main wave},
all the coefficients in the expansion of the right hand side of equation \eqref{eqn:main wave} are non-negative integers,
which is the necessary condition for the existence of counting meaning.

The wave function is closely related to the generating function of superpolynomials of reduced colored HOMFLY-PT homology considered by Sto\v si\'c and Su\l kowski in \cite{SS24} for the torus knot $T_{1,f}$ case.
In their paper,
they proposed an identity between the generating function of standard generalized Schr\"oder paths and a certain ratio of superpolynomials of reduced colored HOMFLY-PT homology of torus knot.
We also give a combinatorial proof of their result for the $T_{1,f}$ case as follows.
\begin{prop}\label{thm:main2}
    Denote by $y_{i}(a,q;x), 1\leq i \leq f+1$ the generating function of generalized Schr\"oder paths under the line $y=\frac{x}{f}$ from $(0,0)$ to $(fl+i-1,l)$ for some integer $l\geq0$, and $\overline{P}^{T_{1,f}}(a,q;x)$ the generating function of superpolynomials of the unreduced colored HOMFLY-PT homology of torus knot $T_{1,f}$.
    There is the following correspondence
    \begin{align}
        y_{i}(a,q;x)
        =\frac{\overline{P}^{T_{1,f}}(a,q;q^{2i-1}x)}
    {\overline{P}^{T_{1,f}}(a,q;q^{-1}x)},
    \quad\quad 1\leq i \leq f+1.
    \end{align}
    See Section \ref{subsec:homology} for more details about the notations.
\end{prop}

The $i=1$ case of above Proposition \ref{thm:main2} is exactly the Proposition 2 in Sto\v si\'c and Su\l kowski's paper \cite{SS24}.
They verified some initial terms for their proposition in Section 4 (see Table 1 in \cite{SS24} for detailed data).

The rest of this paper is organized as follows.
In Section \ref{sec:gSchroder},
we study various kinds of generalized Schr\"oder paths and derive $q$-difference equations for their generating functions,
which proves Theorem \ref{thm:main1}.
In Section \ref{sec:main},
we prove Proposition \ref{thm:main wave},
which establishes the relation between the wave function corresponding to torus knot $T_{1,f}$ and the generating function of generalized Schr\"oder paths with backwards.
In the same section, we also give a combinatorial proof of the Proposition \ref{thm:main2}.

\vspace{.2in}
{\em Acknowledgements}.
The authors would like to thank Sto{\v s}i{\'c} and Su{\l }kowski for their interests in this work.
The first and the third authors thank Jian Zhou for helpful discussions.
The second author is supported by the GRF grant (No. 17303420) of the University Grants Committee of the Hong Kong SAR, China.
The third author is supported by the NSFC grant (No. 12288201),
the China Postdoctoral Science Foundation (No. 2023M743717),
and the China National Postdoctoral Program for Innovative Talents (No. BX20240407).
\vspace{.2in}

 \section{Generalized Schr{\"o}der paths under rational slope}
\label{sec:gSchroder}

In this section,
we study Schr\"oder paths and their various kinds of generalizations.
Then 
we derive the $q$-difference equations that are satisfied by their generating functions.

\subsection{\texorpdfstring{The slope $m/n$}{m/n} case}\label{subsec:slopemn}
A standard generalized Schr{\"o}der path is a lattice path from $(0,0)$ to $(nl,ml)$ with right step $(1,0)$, up step $(0,1)$ and diagonal step $(1,1)$ that never goes above the line $y=\frac{m x}{n}$ for some nonnegative integer $l$. We denote $f=mn$ and $l$ as the size of this path.

For $i,j,l,s\in\mathbb{N}$, we denote $n_{i,j,l}^{[s]}$ as the number of generalized Schr{\"o}der paths that satisfy the following conditions: they have $\frac{i}{2}$ diagonal steps, form an area of $\frac{j}{2}$ between the path and the diagonal line $y=\frac{mx}{n}$, have a size of $l$, and the midway points,
i.e., the points on Schr{\"o}der path where the $x$-coordinate is not equal to $nl$ and the $y$-coordinate is not equal to $0$, 
do not exceed the line $y=\frac{mx-s}{n}$.
For example, if the midway points are empty, then the above conditions are automatically satisfied, therefore for natural number $s\in(mnl-(m+n),mnl]$,
we have
\begin{align}
n^{[s]}_{i, j, l} = 
\left\{
\begin{array}{ll}
     1, & \text{if\ } l \geq 0, i = 0, j = mnl^2,\\
     1, & \text{if\ } l \geq 1, i = 2, j = mnl^2 - 1,\\
     0, & \text{otherwise.}
\end{array}\right.
\end{align}
For $s>mnl$, we adopt the convention that $n_{i,j,l}^{[s]}=0$.
For $(m,n)=(2,3)$ and $l=1$ we have the following examples
\begin{figure}[h!]
\centering
\includegraphics[width=0.9\textwidth]{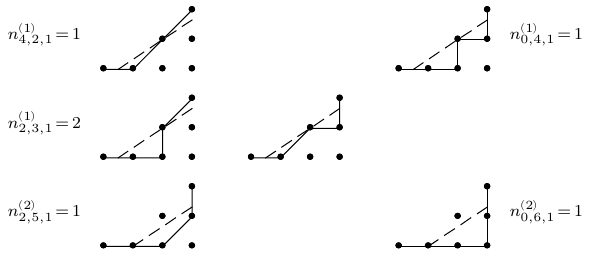}
\caption{Examples for $(m,n)=(2,3)$}
\end{figure}

Let us introduce the following generating function
\begin{align*}
y^{[s]}(a,q;x):=
\sum_{i,j,l\in\mathbb{N}}
n_{i,j,l}^{[s]} a^i q^j x^l.
\end{align*}

\begin{thm}[=Theorem \ref{thm:main1}]
    Denote $\epsilon_s:=\frac{\alpha_sm+\beta_sn-s}{mn}$, where $\alpha_s,\beta_s\in\mathbb{N}$ are determined by
    \begin{alignat*}{2}
    0 &\leq \alpha_s < n,&\quad
    \alpha_sm &= s \mod n,\\
    0 &\leq \beta_s < m,&\quad
    \beta_sn &= s \mod m.
    \end{alignat*}
    The generating functions $y^{[0]},y^{[1]},\ldots,y^{[mn-1]}$ satisfy  the following $q$-difference equations
\begin{subequations}
\begin{align}
    y^{[0]}(a,q;x) & = 
    \frac{1}{1-y^{[1]}(a,q;x)}, \label{eqn:main eqa}\\
    \nonumber
    y^{[s]}(a,q;x) & =
    y^{[s+1]}(a,q;x) +
    q^{mn\epsilon_s^2-2\alpha_s\beta_s}
    x^{-\epsilon_s}\cdot
    y^{[0]}(a,q;q^{2s}x)\\
    & \quad \quad\ \ \cdot
    y^{[\beta_sn+1]}(a,q;q^{2(s-\beta_sn)}x)\cdot
    y^{[\alpha_sm+1]}(a,q;q^{2(s-\alpha_sm)}x), \label{eqn:y^s and yyy}\\
    y^{[mn]}(a,q;x) & = a^2q^{mn-1}x+
    q^{mn}x\cdot y^{[0]}(a,q;q^{2mn}x) \label{eqn:main eqc}
\end{align}
\end{subequations}
for $0<s<mn$.
Moreover,
they are uniquely determined by these $q$-difference equations and the initial values
\begin{align}
    y^{[0]}(a,q;0)=1,\quad
    y^{[s]}(a,q;0)=0,\ 0<s<mn, \quad
    y^{[mn]}(a,q;0)=0.
\end{align}
\end{thm}
{\bf Proof:}
This theorem can be reduced to certain recurrence relation on the numbers $n_{i,j,l}^{[s]}$.
For given $i,j,l,s$ and counting $n_{i,j,l}^{[s]}$,
we consider two cases.
First, we consider the case when the line $y=\frac{mx-s}{n}$ does not pass through any midway points, in which case all the midway points do not exceed the line $y=\frac{mx-(s+1)}{n}$.
Second, we consider the case when the line $y=\frac{mx-s}{n}$ passes through some midway points, in which case any generalized Schr\"oder paths counted by $n_{i,j,l}^{[s]}$ can be uniquely decomposed into the following three parts
\begin{itemize}
    \item[i).] Consider a segment connecting $(\frac{s}{m}, 0)$ with a midway point on the line $y=\frac{mx-s}{n}$, and this segment does not pass through other midway points. If $(\frac{s}{m}, 0)$ itself is a lattice point, then this segment degenerates. Count all paths that the midway points do not exceed this segment, i.e., $n_{i_-,j_-,l_-}^{[\beta_sn+1]}$, where $l_-\geq 1$ and $\beta_s$ is determined by
    \begin{align}
    0 \leq \beta_s < m,\quad
    \beta_sn = s \mod m;
    \end{align}

    \item[ii).] Consider a segment connecting the lattice points on the line $y=\frac{mx-s}{n}$.
    Count all paths that the midway points do not exceed this segment, i.e. $n_{i_0,j_0,l_0}^{[0]}$;

    \item[iii).] Consider a segment connecting a midway point on the line $y=\frac{mx-s}{n}$ with $(nl,ml-\frac{s}{n})$, and this segment does not pass through other midway points on this line. If $(nl,ml-\frac{s}{n})$ itself is a lattice point, then this segment degenerates. Count all paths that the midway points do not exceed this segment, i.e. $n_{i_+,j_+,l_+}^{[\alpha_sm+1]}$, where $l_+\geq 1$ and $\alpha_s$ is determined by
    \begin{align}
    0 \leq \alpha_s < n,\quad
    \alpha_sm = s \mod n.
    \end{align}
\end{itemize}
We can summarize the above decomposition for $n_{i,j,l}^{[s]}$ with given $i,j,l\in\mathbb{N}$ 
and $0\leq s <mn$ as follows
\begin{subequations}
    \begin{align}\nonumber
    n_{i,j,l}^{[s]}  =&
    n_{i,j,l}^{[s+1]} +
    \delta_{i0}\delta_{j0}
    \delta_{l0}\delta_{s0} \\
    \label{nsRec}
    &+\sum_{\substack{
    i_-+i_0+i_+=i\\
    j_-+j_0+j_+=j-(2sl-\frac{s^2}{mn})\\
    l_-+l_0+l_+=l+\epsilon_s\\
    l_-,l_+\geq 1}}
    n_{i_-,j_-+(2\beta_snl_--\frac{\beta_s^2n}{m}),l_-    }^{[\beta_sn+1]}
    n_{i_0,j_0,l_0}^{[0]}
    n_{i_+,j_++(2\alpha_sml_+-\frac{\alpha_s^2m}{n}),l_+}^{[\alpha_sm+1]},\\
    \label{nmnRec}
    n_{i,j,l}^{[mn]}  =&
    n_{i,j-(2l-1)mn,l-1}^{[0]} + 
    \delta_{2i}\delta_{mn-1,j}\delta_{1l}
    .
    \end{align}
\end{subequations}
The above recurrence relations are exactly equivalent to the $q$-difference equations \eqref{eqn:main eqa}, \eqref{eqn:y^s and yyy} and \eqref{eqn:main eqc}.
Moreover,
since $\alpha_sm+\beta_sn-s\in(-mn,2mn)$, the integer $\epsilon_s$ can only be $0$ or $1$.
For $\epsilon_s=0$, we have $l_-,l_0,l_+<l$.
For $\epsilon_s=1$, we can also decompose it into the following three smaller parts
\begin{align*}
    \beta_sn+1>s,\quad l_-\leq l,\quad
    l_0<l,\quad
    \alpha_sm+1>s,\quad l_+\leq l.
\end{align*}
Therefore starting from \eqref{nsRec} and \eqref{nmnRec}, we can always reduce the calculation of $n_{i,j,l}^{[s]}$ to some smaller cases, which ensures that the sequence $(n_{i,j,l}^{[s]})_{i,j,l\in\mathbb{N}}$ is uniquely determined.
Thus, we complete the proof.
$\Box$

Specifically, if we take $m\mid s$ or $n\mid s$ for $0<s<mn$, then equation \eqref{eqn:y^s and yyy} gives
\begin{align}\label{Simpys}
y^{[s]}(a,q;x)= y^{[s+1]}(a,q;x)\cdot y^{[0]}(a,q;q^{2s}x).
\end{align}
Notice that this equation holds for all $s$ when $m=1$.


\subsection{The slope \texorpdfstring{$1/f$}{1/f} case}
\label{sec:gSchroder2}
For $(m, n) = (1, f)$, we introduce the following generating functions for sequences $\{n_{i,j,l;k}\}_{i,j,l\in\mathbb{N},k\in\mathbb{Z}^+}$, $\{n_{i,j,l;\infty}\}_{i,j,l\in\mathbb{N}}$ and $\{u_{i,j,l}\}_{i,j,l\in\mathbb{N}}$
\begin{align*}
y_k(a,q;x) & :=
\sum_{i,j,l\in\mathbb{N}} n_{i,j,l;k} a^iq^jx^l,\\
y_\infty(a,q;x) & :=
\sum_{i,j,l\in\mathbb{N}} n_{i,j,l;\infty} a^iq^jx^l,\\
h(a,q;x) & :=
\sum_{i,j,l\in\mathbb{N}} u_{i,j,l} a^iq^jx^l.
\end{align*}

The integer $n_{i,j,l;k}$ counts the lattice paths that satisfy the following conditions
\begin{itemize}
    \item From $(0, 0)$ to $(fl+k-1, l)$ with right step $(1, 0)$, up step $(0, 1)$ and $\frac{i}{2}$ diagonal step $(1, 1)$;
    \item With infinite numbers of right step $(1, 0)$ after $(fl+k-1, l)$;
    \item Never go above the line $y = \frac{x}{f}$;
    \item Area enclosed by the path, the line  $y = \frac{x}{f}$, and the line $y=l$ is $\frac{j}{2}$.
\end{itemize}

The integer $n_{i,j,l;\infty}$ counts all the paths that are counted by $n_{i,j,l;k}$ for all $k$.
That is to say,
$n_{i,j,l;\infty}$ counts the paths as mentioned above with $\frac{i}{2}$ diagonal steps, encloses $\frac{j}{2}$ area and has a finite height of $l$.

The integer $u_{i,j,l}$ counts the lattice paths that satisfy the following conditions
\begin{itemize}
    \item From $(0, 0)$ to $(fl, l)$ with left step $(-1,0)$, right step $(1, 0)$, up step $(0, 1)$ and $\frac{i}{2}$ diagonal step $(1, 1)$;
    \item Each row will end with some left steps $(-1,0)$ or ends at $(fl, l)$. After the left step, it will not be followed by a right step, and vice versa;
    \item Never go above the line $y = \frac{x}{f}$;
    \item Area enclosed by the path and the line  $y = \frac{x}{f}$ is $\frac{j}{2}$.
\end{itemize}

\begin{figure}[h!]
\centering
\includegraphics[width=1\textwidth]{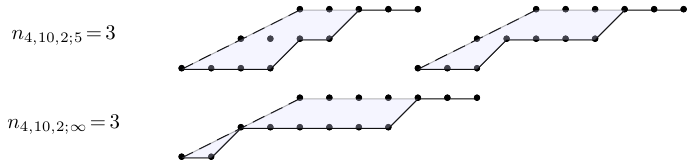}
\caption{Examples for $n_{4,10,2;5}$ and $n_{4,10,2;\infty}$ when $f=2$}
\end{figure}

\begin{figure}[h!]
\centering
\includegraphics[width=0.86\textwidth]{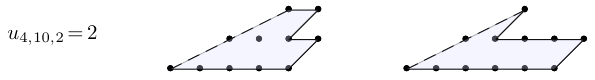}
\caption{Examples for $u_{4,10,2}$ when $f=2$}
\end{figure}
\noindent We call this kind of paths the generalized Schr\"oder paths with backwards since it admits the left step $(-1,0)$.

About the relation between the generating functions $\{y_k(a,q;x)\}_{k\in\mathbb{Z}^{+}}$ and $y_\infty(a,q;x)$,
we have the following lemma.
\begin{lem}\label{lem:yk to yinf}
    For $k\in\mathbb{Z}^+$ we have
    \begin{align}
        n_{i,j,l;k}\leq
        n_{i,j,l;k+1}\leq
        \max\{(j-f+2)^l,0\},
    \end{align}
    and the following limits
    \begin{align}
        n_{i,j,l;\infty}=
        \lim_{k\to\infty}n_{i,j,l;k},\quad
        y_{\infty}(a,q;x)=\lim_{k\to\infty} y_k(a,q;x).
    \end{align}
\end{lem}
{\bf Proof:}
From the above definition, we can see that $n_{i,j,l;k+1}$ also counts all the paths counted by $n_{i,j,l;k}$. Therefore $n_{i,j,l;k+1}\geq n_{i,j,l;k}$.

Since every shape enclosed by a path of size $l$ can be decomposed into $l$ rows, and the area of each row can only be $\frac{f-1}{2},\frac{f}{2},\ldots,\frac{j-1}{2}$, or $\frac{j}{2}$, and each area corresponds to one shape at most, therefore for any $k\in\mathbb{Z}^+$, we have $$n_{i,j,l;k}\leq \max\{(j-f+2)^l,0\}.$$
$\Box$

\begin{prop}
\label{prop:y_1,y_k}
The generating functions $\{y_{k}(a,q;x)\}_{k\in\mathbb{Z}^+}$
are determined by the following $q$-difference equations
\begin{subequations}
\begin{align}
    y_1(a,q;x) & =
    q^fx\cdot y_{f+1}(a,q;x)+
    a^2q^{f-1}x\cdot y_{f}(a,q;x)+1,
    \label{eqn:y_1 q-diff}\\
    y_k(a,q;x)&=y_1(a,q;q^{2(k-1)}x)\cdot 
    y_{k-1}(a,q;x),\quad
    k>1, \label{eqn:y_k q-diff}
\end{align}
with the initial value $y_1(a,q;0)=1$.
\end{subequations}
The generating function $y_{\infty}(a,q;x)$ is determined by the following $q$-difference equation
\begin{align}\label{eqn:y_inf q-diff}
\frac{q^fx}{y_\infty(a,q;q^{2(f+1)}x)} +
\frac{a^2q^{f-1}x}{y_\infty(a,q;q^{2f}x)} -
\frac{1}{y_\infty(a,q;q^{2}x)} +
\frac{1}{y_\infty(a,q;x)} = 0,
\end{align}
with the initial value $y_\infty(a,q;0)=1$.
\end{prop}
{\bf Proof:}
First, notice that we can divide the paths counted by $n_{i,j,l+1;1}$ into the following two classes.
\begin{itemize}
    \item[i).] Ends with an up step $(0,1)$ starting from $(fl+f,l)$, which has $n_{i,j-f,l;f+1}$ paths;
    \item[ii).] Ends with a diagonal step $(1,1)$ starting from $(fl+f-1,l)$, which has $n_{i-2,j-(f-1),l;f}$ paths.
\end{itemize}
Thus,
for $l\in\mathbb{N}$ we have
\begin{align*}
n_{i,j,l+1;1}=
n_{i,j-f,l;f+1} +
n_{i-2,j-(f-1),l;f}.
\end{align*}
In addition, there is only one path of size $0$.
Therefore, we have completed the proof of equation \eqref{eqn:y_1 q-diff}.

Next, notice that we can uniquely decompose a path from $(0,0)$ to $(fl_{k}+k-1,l_{k})$, counted by $n_{i,j,l_k;k}$, into the following parts.
\begin{itemize}
    \item[i).] From $(0,0)$ to $(fl_{k-1}+k-2,l_{k-1})$, where $l_{k-1}\leq l_k$;
    \item[ii).] From $(fl_{k-1}+k-2,l_{k-1})$ to $(fl_{k-1}+k-1,l_{k-1})$;
    \item[iii).] From $(fl_{k-1}+k-1,l_{k-1})$ to $(fl_{k}+k-1,l_{k})$, and this part does not exceed the line $y=\frac{x-k+1}{f}$.
    In this way, $l_{k-1}$ becomes unique.
\end{itemize}
Notice that the paths in the third part can be obtained by shifting the paths counted by $y_1(a,q;x)$ to the right by $k-1$ units, thus the contribution from this part is the generating series $y_1(a,q;q^{2(k-1)}x)$.
Combining all the above discussions,
we have completed the proof of equation \eqref{eqn:y_k q-diff}. 

As a corollary of equation \eqref{eqn:y_k q-diff} and Lemma \ref{lem:yk to yinf}, we have
\begin{align*}
    y_k(a,q;x) =
    \frac{y_\infty(a,q;x)}
    {y_\infty(a,q;q^{2k}x)}.
\end{align*}
Thus, equation \eqref{eqn:y_inf q-diff} can be reduced to the equation \eqref{eqn:y_1 q-diff}.
$\Box$


\begin{rmk}
    From the definitions,
    $y_1(a,q;x)$ and $y^{[0]}(a,q;x)$ count the same generalized Schr\"oder paths.
    Thus,
    $y_1(a,q;x)$ can also be calculated from Theorem \ref{thm:main1}
    which leads the same formula in Proposition \ref{prop:y_1,y_k}.
\end{rmk}

\begin{prop}
\label{prop:h}
The generating functions $h(a,q;x)$ of generalized Schr\"oder paths with backwards
are determined by the following $q$-difference equations
\begin{align}\label{eqn:h q-diff}
h(a,q;x) =
h(a,q;q^{2}x) +
q^fx\cdot h(a,q;q^{2(f+1)}x) +
a^2q^{f-1}x\cdot h(a,q;q^{2f}x)
\end{align}
with the initial value $h(a,q;0)=1$.
\end{prop}
{\bf Proof:}
First, notice that we can divide the paths counted by $u_{i,j,l+1}$ into the following three classes:
\begin{itemize}
    \item[i).] Passing through $(f(l+1)+1,l+1)$;
    \item[ii).] From $(fl+f+1,l)$ reach the end point $(f(l+1),l+1)$ through left step $(-1,0)$ and up step $(0,1)$;
    \item[iii).] From $(fl+f,l)$ reach the end point $(f(l+1),l+1)$ through left step $(-1,0)$ and diagonal step $(1,1)$.
\end{itemize}
The paths of the first part can uniquely come from the paths with backwards moving $1$ unit to the right, so the contribution from this part is $h(a,q;q^2x)$;
The paths of the second part can uniquely come from the paths with backwards moving $f+1$ units to the right,
so the contribution from this part is $q^fx\cdot h(a,q;q^{2(f+1)}x)$;
The paths of the third part can uniquely come from the counted path moving $f$ units to the right,
so the contribution from this part is $a^2q^{f-1}x\cdot h(a,q;q^{2f}x)$.
Thus, we complete the proof.
$\Box$

\begin{cor}
    The generating functions $y_\infty(a,q;x)$ and $h(a,q;x)$ are related to each other by the following
    $$h(a,q;x)=\frac{1}{y_\infty(a,q;-x)}.$$
\end{cor}
{\bf Proof:}
Comparing equations \eqref{prop:y_1,y_k} and \eqref{eqn:h q-diff}, one can see that $h(a,q;x)$ and $\frac{1}{y_\infty(a,q;-x)}$ satisfy the same $q$-difference equation and initial value, thus they are equal.
$\Box$

\section{Knot invariants and counting generalized Schr{\"o}der paths}
\label{sec:main}

\subsection{Wave function corresponding to $T_{1,f}$ and generalized Schr{\"o}der paths with backwards}\label{subsec:wave}
In this subsection,
we establish the relation between the wave function of the partition function for colored HOMFLY-PT polynomials of torus knot $T_{1,f}$ and the generating function of generalized Schr{\"o}der paths with backwards,
which proves Proposition \ref{thm:main wave}.

The colored HOMFLY-PT polynomials for torus knot $T_{m,n}$ were widely studied.
They were explicitly computed by Rosso and Jones \cite{RJ93} using the language of the quantum group.
Later,
Brini, Eynard and Mari\~no \cite{BEM12} derived the spectral curve under the large $N$ duality and mirror symmetry, see also \cite{FZ19}.
Moreover,
the KP integrability of the partition function for these knot invariants was studied in \cite{MMS13, DPSS19}.
In this paper,
we follow the notations in \cite{LZ10},
in where the authors provide explicit formulas for these colored HOMFLY-PT polynomials of torus knots in terms of Schur polynomials $s_\lambda(\mathbf{t})$ and related combinatorics.
We briefly review the notion of Schur polynomial in Appendix \ref{sec:app} for the reader's convenience.

Following the notation in \cite{LZ10},
for two relatively prime integers $(m,n)$,
the colored HOMFLY-PT polynomial of torus knot $T_{m,n}$ labeled by a partition $\lambda$ is given by (here we consider the torus knot $T_{m,n}$ with framing $mn$)
\begin{align}\label{eqn:H_lambda}
    H_\lambda^{T_{m,n}}(\nu,t)
    =\nu^{n(m-1)|\lambda|/2}
    \sum_{\mu\vdash m|\lambda|}
    C^{\lambda}_{\mu,m} t^{-\kappa_{\mu}n/2m}
    s_\mu(\mathbf{t}^*),
\end{align}
where
$\mu\vdash m|\lambda|$ means
$\mu$ is partition of $m|\lambda|$,
$\kappa_{\mu}:=2\sum_{i=1}^{l(\mu)}\sum_{j=1}^{\mu_i} (j-i)$,
and $s_\mu(\mathbf{t}^*)$ denotes the evaluation of Schur polynomial at $t_i=t^*_i:=\frac{1}{i}\cdot\frac{\nu^{i/2}-\nu^{-i/2}}{t^{i/2}-t^{-i/2}}$.
The coefficients $C^{\lambda}_{\mu,m}$ are determined by
\begin{align*}
    s_{\lambda}(mt_{m},mt_{2m},mt_{3m}...)
    =\sum_{\mu\vdash m|\lambda|}
    C^{\lambda}_{\mu,m} \cdot s_{\mu}(t_1,t_2,t_3,...).
\end{align*}

The partition function corresponding to these colored HOMFLY-PT polynomials is defined by
\begin{align*}
    Z^{T_{m,n}}(\mathbf{t})
    :=\sum_{\lambda\in\mathcal{P}}
    H_\lambda^{T_{m,n}}(\nu,t) s_{\lambda}(\mathbf{t}).
\end{align*}
The wave function corresponding to this partition function $Z^{T_{m,n}}(\mathbf{t})$ is defined by
(see, for example, Section 6 in \cite{DJM00})
\begin{align*}
    \psi^{T_{m,n}}(x):=\frac{Z^{T_{m,n}}(\mathbf{t})|_{t_i=x^i/i}}{Z^{T_{m,n}}(0)}
    =\sum_{k\geq0}
    H_{(k)}^{T_{m,n}}(\nu,t) x^k,
\end{align*}
where we have used that $Z^{T_{m,n}}(0)=1$ and
$s_{\lambda}(\mathbf{t})|_{t_i=x^i/i}=\delta_{\lambda,(k)} x^k$
for some $k\geq0$.
As a consequence,
the wave function $\psi^{T_{m,n}}(x)$ is the generating function of the colored HOMFLY-PT polynomials of torus knot labeled by one-row partitions.

\begin{prop}[=Proposition \ref{thm:main wave}]
    Denote by $\psi^{T_{1,f}}(x)$ the wave function corresponding to the partition function of colored HOMFLY-PT polynomials of torus knot $T_{1,f}$,
    we have
    \begin{align}
        h(a,q;x)
        =\psi^{T_{1,f}}(-\nu^{1/2}q^{f-1}x)|_{\nu\rightarrow-a^2/q, t\rightarrow q^{-2}}.
    \end{align}
\end{prop}
{\bf Proof:}
When restricting to the $T_{1,f}$ case,
these colored HOMFLY-PT polynomials and partition function $Z^{T_{m,n}}(\mathbf{t})$ become much clearer. 
From equation \eqref{eqn:H_lambda},
we have
\begin{align*}
    Z^{T_{1,f}}(\mathbf{t})
    =\sum_{\lambda\in\mathcal{P}}
    t^{-f\kappa_{\lambda}/2}s_{\lambda}(\mathbf{t}^*) s_{\lambda}(\mathbf{t}).
\end{align*}
The wave function related to this partition function is thus
\begin{align*}
    \psi^{T_{1,f}}(x)
    =&Z^{T_{1,f}}(\mathbf{t})|_{t_i=x^i/i}
    =\sum_{k\geq0}
    t^{-f(k^2-k)/2}s_{(k)}(\mathbf{t}^*) \cdot x^k\\
    =&t^{-f\big((x\partial_x)^2-x\partial_x\big)/2}
    \exp\big(\sum_{k\geq1} t_k^* x^k\big),
\end{align*}
where in the last line,
we have used the Cauchy identity for Schur polynomials.
For the exponential function part in the above equation
\begin{align*}
    A(x):=
    \exp\Big(\sum_{k\geq1} t_k^* x^k\Big)
    =\exp\Bigg(\sum_{k\geq1} \frac{(\nu^{k/2}-\nu^{-k/2})x^k}{(t^{k/2}-t^{-k/2})k}\Bigg),
\end{align*}
one can directly verify that it satisfies the following equation
\begin{align*}
    \nu^{1/2}t^{1/2}x A(tx)
    -\nu^{-1/2}t^{1/2}x A(x)
    -A(tx)
    +A(x) = 0.
\end{align*}
As a consequence,
the wave function $\psi(x)$ satisfies the following $q$-difference equation
\begin{align}
    \nu^{1/2}t^{1/2}x \psi^{T_{1,f}}(t^{-f+1}x)
    -\nu^{-1/2}t^{1/2}x \psi^{T_{1,f}}(t^{-f}x)
    -\psi^{T_{1,f}}(tx)
    +\psi^{T_{1,f}}(x) = 0.
\end{align}
Substituting the $\psi^{T_{1,f}}(x)$ in above $q$-difference equation by
$$\psi^{T_{1,f}}(-\nu^{1/2}q^{f-1}x)|_{\nu\rightarrow-a^2/q, t\rightarrow q^{-2}},$$
the result is exactly equal to the equation \eqref{eqn:h q-diff} satisfied by $h(a,q;x)$.
We thus prove this proposition.
$\Box$

\begin{rmk}
Actually,
the change of coordinates used in the Proposition \ref{thm:main wave} makes all the coefficients in the expansion of $\psi^{T_{1,f}}(-\nu^{1/2}q^{f-1}x)|_{\nu\rightarrow-a^2/q, t\rightarrow q^{-2}}$ are non-negative integers,
which is somehow a necessary condition for admitting a counting interpretation.
\end{rmk}

\subsection{Reduced colored HOMFLY-PT homology of $T_{1,f}$ and generalized Schr{\"o}der paths}
\label{subsec:homology}
In this subsection,
we prove Proposition \ref{thm:main2},
which establishes the relation,
proposed by Sto\v si\'c and Su\l kowski \cite{SS24},
between the superpolynomials of reduced colored HOMFLY-PT homology of torus knot $T_{1,f}$ and the generating function of generalized Schr{\"o}der paths.

We directly follow the notations in \cite{SS24}.
The generating function of superpolynomials of unreduced colored HOMFLY-PT homology of the torus knot $T_{1,f}$ is defined by (see equation (29) in \cite{SS24})
\begin{align*}
    \overline{P}^{T_{1,f}}(a,q;x)=&
    \sum_{r=0}^{\infty}
    \overline{P}^{T_{1,f}}_r(a,q) x^r\\
    :=&1+\sum_{r=1}^{\infty}
    (-1)^r\frac{\prod_{i=1}^r(a^2+q^{2i-1})}{\prod_{i=1}^r(1-q^{2i})} q^{fr^2}
    \mathcal{P}_{r}^{T_{1,f}}(a,1,q^{-1}) x^r,
\end{align*}
where $\mathcal{P}^{T_{1,f}}_r(a,q;t), r\geq1$ are the superpolynomials of the reduced colored HOMFLY-PT homology
(see, for examples, \cite{GGS18,SS24}).

\begin{prop}[=Proposition \ref{thm:main2}]
    Denote by $\overline{P}^{T_{1,f}}(a,q;x)$ the generating function of superpolynomials of unreduced colored HOMFLY-PT homology of torus knot $T_{1,f}$,
    we have
    \begin{align}
        y_i(a,q;x)
        =\frac{\overline{P}^{T_{1,f}}(a,q;q^{2i-1}x)}
    {\overline{P}^{T_{1,f}}(a,q;q^{-1}x)},
    1\leq i \leq f+1.
    \end{align}
\end{prop}
{\bf Proof:}
Denote
\begin{align}
    \tilde{y}_i(a,q;x)
    :=\frac{\overline{P}^{T_{1,f}}(a,q;q^{2i-1}x)}
    {\overline{P}^{T_{1,f}}(a,q;q^{-1}x)}
\end{align}
for $1\leq i \leq f+1$.
It is apparent that $\{\tilde{y}_i(a,q;x)\}_{2\leq i\leq f+1}$ can be determined from the following recursive relations
\begin{align}\label{eqn:ty_i q-diff}
    \tilde{y}_i(a,q;x)
    =q^{2(i-1)x\partial_x}\tilde{y}_1(a,q;x) \cdot \tilde{y}_{i-1}(a,q;x),
    \quad\quad 2\leq i\leq f+1,
\end{align}
which is exactly equivalent to those recursive equations \eqref{eqn:y_k q-diff} satisfied by $y_i(a,q;x), 2\leq i\leq f+1$.
Thus,
We only need to show that $\tilde{y}_1(a,q;x)$ is equal to $y_1(a,q;x)$.
We will prove this fact by comparing the $q$-difference equations that they satisfy.

For the torus knot $T_{1,f}$,
it is well-known that the superpolynomials of unreduced colored HOMFLY-PT homology are determined by (see equation (73) in \cite{SS24})
\begin{align}\label{eqn:rec for p_r}
    \overline{P}_{r+1}^{T_{1,f}}(a,q)
    =-q^{(2r+1)f} \cdot \frac{a^2+q^{2r+1}}{1-q^{2r+2}}
    \cdot \overline{P}_{r}^{T_{1,f}}(a,q),
    \quad\quad r\geq0
\end{align}
and $\overline{P}_{0}^{T_{1,f}}(a,q)=1$.
This recursive relation \eqref{eqn:rec for p_r} is equivalent to the following $q$-difference equation for the generating function $\overline{P}^{T_{1,f}}(a,q;x)$,
\begin{align}\label{eqn:Pbar q-diff}
\begin{split}
    \overline{P}^{T_{1,f}}&(a,q;q^{2f+2}x)\\
    =&\Big(\overline{P}^{T_{1,f}}(a,q;q^2x)
    -q^{f} a^2 x \overline{P}^{T_{1,f}}(a,q;q^{2f}x)
    -\overline{P}^{T_{1,f}}(a,q;x)\Big)/(q^{f+1} x).
\end{split}
\end{align}
The above equation is equivalent to the equation (75) in \cite{SS24} regarding quantum A-polynomial.
Consider
\begin{align}
    &q^{(2f-2i+2)x\partial_x} \cdot \tilde{y}_i(a,q;x) \nonumber\\
    =&\frac{\overline{P}^{T_{1,f}}(a,q;qx)
    -q^{f-1} a^2 x \overline{P}^{T_{1,f}}(a,q;q^{2f-1}x)
    -\overline{P}^{T_{1,f}}(a,q;q^{-1}x)}
    {q^f x \overline{P}^{T_{1,f}}(a,q;q^{2f-2i+1}x)}\nonumber \\
    =&\frac{1}{q^f x}
    \Big(q^{2x\partial_x} \frac{1}{\tilde{y}_{f-i}(a,q;x)}
    -q^{f-1} a^2 x q^{(2f-2i+2)x\partial_x} \tilde{y}_{i-1}(a,q;x)
    -\frac{1}{\tilde{y}_{f-i+1}(a,q;x)}
    \Big),
\end{align}
where $\tilde{y}_{0}(a,q;x):=1$ and in the first equal sign,
we have used the $q$-difference equation \eqref{eqn:Pbar q-diff}.
By setting $i=f$,
the above equation gives the following $q$-difference equation
\begin{align}\label{eqn:last eqn}
    q^f x \cdot q^{2x\partial_x} \cdot \tilde{y}_f(a,q;x)
    =&1
    -q^{f-1} a^2 x \cdot q^{2x\partial_x} \tilde{y}_{f-1}(a,q;x)
    -\frac{1}{\tilde{y}_{1}(a,q;x)}.
\end{align}
By multiplying $\tilde{y}_{1}(a,q;x)$ to both sides of above equation and using equation \eqref{eqn:ty_i q-diff},
the above equation \eqref{eqn:last eqn} is exactly equal to the $q$-difference equation \eqref{eqn:y_1 q-diff} satisfied by $\{y_i(a,q;x)\}_{1\leq i\leq f+1}$.
Thus,
this proposition is proved.
$\Box$


\appendix

\section{Schur polynomials}
\label{sec:app}
In this appendix,
we review the notion of Schur polynomials for convenience and completeness.
We recommend the book \cite{Mac15}.

A partition $\lambda=(\lambda_1,...,\lambda_{l(\lambda)})$ is a series of non-increasing positive integers.
The length and size of the partition $\lambda$ are denoted by $l(\lambda)$ and $|\lambda|:=\sum_{i=1}^{l(\lambda)} \lambda_i$, respectively.
The elementary Schur polynomials $s_{(n)}(\mathbf{t})\in\mathbb{C}[t_1,t_2,...]$ are defined,
in terms of their generating function,
by
\begin{align*}
    \sum_{n=0}^\infty
    s_{(n)}(\mathbf{t}) x^n
    :=\exp\big(\sum_{k=1}^\infty t_k x^k \big).
\end{align*}
Then, in general,
the Schur polynomial labeled by a partition $\lambda$ can be obtained by the following Jacobi-Trudi formula
\begin{align*}
    s_{\lambda}(\mathbf{t})
    =\det\big( s_{(\lambda_i-i+j)}(\mathbf{t}) \big)_{1\leq i,j\leq l(\lambda)}\in\mathbb{C}[t_1,t_2,...].
\end{align*}

Certain special evaluations of Schur polynomials are important in the areas of representation theory and combinatorics
(see Section I in \cite{Mac15} for elaborate results about special evaluations of Schur polynomials).
In this paper,
we only need the following
\begin{align*}
    s_{\lambda}(\mathbf{t})|_{t_k=x^k/k}
    =x^n
\end{align*}
if $\lambda=(n)$ is a one-row partition for some non-negative integer $n$,
otherwise it is zero.
This evaluation is related to the wave function and quantum spectral curve (see \cite{DJM00} and \cite{GS12}).

\renewcommand{\refname}{Reference}
\bibliographystyle{plain}
\bibliography{reference}

\begin{thebibliography}{10}

\bibitem{AKMV}
Mina Aganagic, Albrecht Klemm, Marino Mari\~{n}o, and Cumrun Vafa.
\newblock The topological vertex.
\newblock {\em Comm. Math. Phys.}, 254(2):425--478, 2005.

\bibitem{BF02}
Cyril Banderier and Philippe Flajolet.
\newblock Basic analytic combinatorics of directed lattice paths.
\newblock volume 281, pages 37--80. 2002.
\newblock Selected papers in honour of Maurice Nivat.

\bibitem{BW19}
Cyril Banderier and Michael Wallner.
\newblock The kernel method for lattice paths below a line of rational slope.
\newblock In {\em Lattice path combinatorics and applications}, volume~58 of
  {\em Dev. Math.}, pages 119--154. Springer, Cham, 2019.

\bibitem{Biz54}
Michael T.~L. Bizley.
\newblock Derivation of a new formula for the number of minimal lattice paths
  from {$(0,0)$} to {$(km,kn)$} having just {$t$} contacts with the line
  {$my=nx$} and having no points above this line; and a proof of {G}rossman's
  formula for the number of paths which may touch but do not rise above this
  line.
\newblock {\em J. Inst. Actuar.}, 80:55--62, 1954.

\bibitem{BEM12}
Andrea Brini, Bertrand Eynard, and Marcos Mari\~no.
\newblock Torus knots and mirror symmetry.
\newblock {\em Ann. Henri Poincar\'e}, 13(8):1873--1910, 2012.

\bibitem{DPSS19}
Petr Dunin-Barkowski, Aleksandr Popolitov, Sergey Shadrin, and Alexey Sleptsov.
\newblock Combinatorial structure of colored {HOMFLY}-{PT} polynomials for
  torus knots.
\newblock {\em Commun. Number Theory Phys.}, 13(4):763--826, 2019.

\bibitem{FZ19}
Bohan Fang and Zhengyu Zong.
\newblock Topological recursion for the conifold transition of a torus knot.
\newblock {\em Selecta Math. (N.S.)}, 25(3):Paper No. 35, 44, 2019.

\bibitem{HOMFLY}
Peter Freyd, David Yetter, Jim Hoste, William B.~R. Lickorish, Kenneth Millett,
  and Adrian Ocneanu.
\newblock A new polynomial invariant of knots and links.
\newblock {\em Bull. Amer. Math. Soc. (N.S.)}, 12(2):239--246, 1985.

\bibitem{GGS18}
Eugene Gorsky, Sergei Gukov, and Marko Sto\v{s}i\'c.
\newblock Quadruply-graded colored homology of knots.
\newblock {\em Fund. Math.}, 243(3):209--299, 2018.

\bibitem{GS12}
Sergei Gukov and Piotr Su{\l }kowski.
\newblock A-polynomial, {B}-model, and quantization.
\newblock {\em J. High Energy Phys.}, (2):070, front matter+56, 2012.

\bibitem{Jon87}
Vaughan F.~R. Jones.
\newblock Hecke algebra representations of braid groups and link polynomials.
\newblock {\em Ann. of Math. (2)}, 126(2):335--388, 1987.

\bibitem{LZ10}
Xiao-Song Lin and Hao Zheng.
\newblock On the {H}ecke algebras and the colored {HOMFLY} polynomial.
\newblock {\em Trans. Amer. Math. Soc.}, 362(1):1--18, 2010.

\bibitem{Mac15}
Ian~G. Macdonald.
\newblock {\em Symmetric functions and {H}all polynomials}.
\newblock Oxford Classic Texts in the Physical Sciences. The Clarendon Press,
  Oxford University Press, New York, second edition, 2015.
\newblock With contribution by A. V. Zelevinsky and a foreword by Richard
  Stanley.

\bibitem{MMS13}
Andrei~D. Mironov, Alexei~Yu. Morozov, and Alexey~V. Sleptsov.
\newblock Genus expansion of {HOMFLY} polynomials.
\newblock {\em Theoret. and Math. Phys.}, 177(2):1435--1470, 2013.
\newblock Russian version appears in Teoret. Mat. Fiz. 177 (2013), no. 2,
  179--221.

\bibitem{DJM00}
Tetsuji Miwa, Michio Jimbo, and Etsuro Date.
\newblock Differential equations, symmetries and infinite dimensional algebras.
\newblock {\em Cambridge University Press}, 2000.

\bibitem{PSS18}
Mi{\l }osz Panfil, Marko Sto{\v s}i{\'c}, and Piotr Su{\l }kowski.
\newblock Donaldson-{T}homas invariants, torus knots, and lattice paths.
\newblock {\em Phys. Rev. D}, 98(2):026022, 28, 2018.

\bibitem{PT}
J\'ozef~H. Przytycki and Pawe\l\ Traczyk.
\newblock Conway algebras and skein equivalence of links.
\newblock {\em Proc. Amer. Math. Soc.}, 100(4):744--748, 1987.

\bibitem{RT90}
Nikolai~Yu. Reshetikhin and Vladimir~G. Turaev.
\newblock Ribbon graphs and their invariants derived from quantum groups.
\newblock {\em Comm. Math. Phys.}, 127(1):1--26, 1990.

\bibitem{RJ93}
Marc Rosso and Vaughan Jones.
\newblock On the invariants of torus knots derived from quantum groups.
\newblock {\em J. Knot Theory Ramifications}, 2(1):97--112, 1993.

\bibitem{Son05}
Chunwei Song.
\newblock The generalized {S}chr\"oder theory.
\newblock {\em Electron. J. Combin.}, 12:Research Paper 53, 10, 2005.

\bibitem{SS24}
Marko Sto{\v s}i{\'c} and Piotr Su{\l }kowski.
\newblock Torus knots and generalized schr\"oder paths.
\newblock {\em arXiv:2405.10161}, 2024.

\bibitem{Wit89}
Edward Witten.
\newblock Quantum field theory and the {J}ones polynomial.
\newblock In {\em Braid group, knot theory and statistical mechanics}, volume~9
  of {\em Adv. Ser. Math. Phys.}, pages 239--329. World Sci. Publ., Teaneck,
  NJ, 1989.

\bibitem{Wit95}
Edward Witten.
\newblock Chern-{S}imons gauge theory as a string theory.
\newblock In {\em The {F}loer memorial volume}, volume 133 of {\em Progr.
  Math.}, pages 637--678. Birkh\"auser, Basel, 1995.

\end{thebibliography}
\vspace{30pt} \noindent
\end{document}